\documentclass[12pt]{article}

\def\edo{\end{document}	 }

\usepackage{amsmath, amsthm, amscd, amsfonts, amssymb}

\newtheorem{theorem}{Theorem}[section]

\newtheorem{lemma}[theorem]{Lemma}
\newtheorem{definition}[theorem]{Definition}

\newtheorem{remark}[theorem]{Remark}

\usepackage[usenames,dvipsnames]{color}

\def\rrd{{\mathbb{R}^d}}
\def\call{{\mathcal{L}}}
\def\calz{{\mathcal{Z}}}

\def\calf{{\mathcal{F}}}

\def\calo{{\mathcal{O}}}

\def\cald{{\mathcal{D}}}
\def\calx{{\mathcal{X}}}

\def\call{{\mathcal{L}}}

\def\vsp{\vspace*{1,5mm}\\ }
\def\bk{\bigskip }
\def\mk{\medskip }
\def\sk{\smallskip }
\def\n{\noindent }
\def\dd{\displaystyle}

\def\barr{\begin{array}}
\def\earr{\end{array}}

\def\bit{\begin{itemize}}
\def\eit{\end{itemize}}

\def\FP{Fokker--Planck}

\def\1{^{-1}}
\def\one{\mbox{$1\!\!\,\rule{0,2mm}{3,1mm}\,$}}

\def\rr{{\mathbb{R}}}
\def\nn{{\mathbb{N}}}
\def\9{{\infty}}
\def\lbb{{\lambda}}
\def\wt{\widetilde}
\def\ov{\overline}
\def\vf{{\varphi}}

\def\ooo{{\Omega}}

\def\vp{{\varepsilon}}

\def\ff{\forall }
\def\({\left(}
\def\){\right)}
\def\<{\left<}
\def\>{\right>}

\title{Nonlinear Fokker--Planck equations  with~fractional Laplacian\\ and McKean--Vlasov SDEs with L\'evy--Noise}
\author{Viorel Barbu\thanks{Al.I. Cuza University and Octav Mayer Institute of Mathematics of  Romanian Academy, Ia\c si, Romania.  Email: vbarbu41@gmail.com}\and Michael R\"ockner\thanks{Fakult\"at f\"ur Mathematik, Universit\"at Bielefeld,  D-33501 Bielefeld, Germany.  Email: roeckner@math.uni-bielefeld.de}}
\date{}

\begin{document}
\maketitle
\begin{abstract}
\n This work is concerned with the existence  of mild solutions to non\-linear Fokker--Planck equations with fractional Laplace operator  $(-\Delta)^s$ for  $s\in\left(\frac12,1\right)$. The uniqueness of Schwartz  distributional solutions is also proved  under suitable assumptions on diffusion and drift terms. As applications, weak existence and uniqueness of solutions to McKean--Vlasov equations with L\'evy--Noise, as well as the Markov property for their laws are proved.\sk\\
{\bf MSC:} 60H15, 47H05, 47J05.\\
{\bf Keywords:} Fokker--Planck  equation, fractional Laplace operator, distributional solutions, mild solution, stochastic differential equation, superposition principle, L\'evy processes.  
\end{abstract}

\section{Introduction}\label{s1}
We  consider here the nonlinear \FP\ equation (NFPE) 
\begin{equation}\label{e1.1}
\barr{l}
u_t+(-\Delta)^s\beta(u)+{\rm div}(Db(u)u)=0,\ \mbox{ in  } (0,\9)\times\rr^d,\\
u(0,x)=u_0(x),\ x\in\rr^d,
\earr\end{equation}where $\beta:\rr\to\rr$, $D:\rrd\to\rrd$, $d\ge2$, and $b:\rr\to\rr$ are given functions to be made precise later on, while $(-\Delta)^s$, $0<s<1$, is the fractional Laplace operator defined as follows. Let $S':=S'(\rrd)$ be the dual of the Schwartz test function space $S:=S(\rrd)$. Define
$$D_s:=\{u\in S';\ \calf(u)\in L^1_{\rm loc}(\rrd),\ |\xi|^{2s}\calf(u)\in S'\}\ (\supset L^1(\rrd))$$ and 
\begin{equation}\label{e1.2}
\calf((-\Delta)^su)(\xi)=|\xi|^{2s}\calf(u)(\xi),\ \xi\in\rrd,\ u\in D_s,\end{equation}where $\calf$ stands for the Fourier transform in $\rrd$, that is,
\begin{equation}\label{e1.3}
\calf(u)(\xi)=(2\pi)^{-d/2}\int_\rrd e^{ix\cdot\xi}u(x)dx,\ \xi\in\rrd,\ u\in L^1(\rrd).\end{equation}($\calf$ extends from $S'$ to itself.)

NFPE \eqref{e1.1} is used for modelling the dynamics of anomalous diffusion of particles in disordered media. The solution $u$ may be viewed as the transition density corresponding to a distribution dependent stochastic differential equation with L\'evy forcing term.

\bk\n{\bf Hypotheses}
\begin{itemize}
	\item[\rm(i)] $\beta\in C^1(\rr)\cap {\rm Lip}(\rr),\ \beta'(r)>0,\ \ff\,r\ne0.$ 
	\item[\rm(ii)] $D\in L^\9(\rrd;\rrd),\ {\rm div}\,D\in L^2_{\rm loc}(\rrd).$ 
	\item[\rm(iii)] $b\in C_b(\rr).$
	\item[\rm(iv)] $({\rm div}\,D)^-\in L^\9,\ b\ge0.$
\end{itemize}	
	Here, we shall study the existence of a {\it mild solution} to equation \eqref{e1.1}\break (see Definition \ref{d1.1} below) and also the uniqueness of distributional solutions. As regards the existence, we shall follow the semigroup methods used in \cite{2}--\cite{4} in the special case $s=1$. Namely, we shall represent  \eqref{e1.1} as an abstract differential equation in $L^1(\rrd)$ of the form
\begin{equation}\label{e1.4}
		\barr{l}
		\dd\frac{du}{dt}+A(u)=0,\ \ t\ge0,\vsp
		u(0)=u_0,\earr \end{equation}where $A$ is a suitable realization in $L^1(\rrd)$ of the operator
\begin{equation}\label{e1.5}
\barr{rl}
A_0(u)\!\!\!&=(-\Delta)^s\beta(u)+{\rm div}(Db(u)u),\  u\in D(A_0),\vsp
D(A_0)\!\!\!&=\left\{u\in L^1(\rrd);(-\Delta)^s\beta(u)+{\rm div}(Db(u)u)\in L^1(\rrd)\right\},\earr\end{equation}	
where div is taken in the sense of Schwartz distributions on $\rrd$.	
	
\begin{definition}\label{d1.1}\rm A function $u\in C([0,\9);L^1:=L^1(\rrd))$ is said to be a {\it mild solution} to \eqref{e1.1} if, for each $0<T<\9$,  
	\begin{eqnarray}
&u(t)=\dd\lim_{h\to0}u_h(t)\mbox{ in }L^1(\rrd),\ t\in[0,T),\label{e1.6}\end{eqnarray}where
\begin{eqnarray}
	&u_h(t)=u^j_h,\ \ff\,t\in [jh,(j+1)h),\ j=0,1,...,N=\mbox{$\left[\frac Th\right]$},\label{e1.7}\\[1mm]
	&u^{j+1}_h+hA_0(u^{j+1}_h)=u^{j}_h,\ j=0,1,....,N_h,\label{e1.8}\\[1mm]
	&u^j_h\in D(A_0),\ \ff\,j=0,...,N_h;\ u^0_h=u_0.\label{e1.9}\end{eqnarray}
\end{definition}	

Of course, Definition \ref{d1.1} makes sense only if the range $R(I+hA_0)$ of the operator $I+hA_0$ is all of $L^1(\rrd)$. 
We note that, if $u$ is a mild solution to \eqref{e1.1}, then it is also a {\it Schwartz distributional solution}, that is,
\begin{equation}\label{e1.10}
\hspace*{-5mm}\barr{l}
\dd\int^\9_0\!\!\!\int_\rrd\!(u(t,x)
\vf_t(t,x)-(-\Delta)^s\vf(t,x)\beta(u(t,x))\\
\qquad\qquad+b(u(t,x))u(t,x)D(x)\cdot\nabla\vf(t,x))dtdx\vsp
\qquad\qquad+\dd\int_\rrd \vf(0,x)u_0(dx)=0,\ \ff\,\vf\in C^\9_0([0,\9)\times\rrd), \earr\end{equation} where  $u_0$ is a measure of finite variation on $\rrd$. 
The main existence result for equation \eqref{e1.1} is given by Theorem \ref{t2.3} below,  which amounts to saying  that under Hypotheses (i)--(iv) there is a mild solution $u$ represented as $u(t)=S(t)u_0,$ $t\ge0$, where $S(t)$ is a continuous semigroup of nonlinear contractions in $L^1$. In Section \ref{s3}, the uniqueness of dis\-tri\-bu\-tional solutions to \eqref{e1.1}, \eqref{e1.10} respectively, in the class $(L^1\cap L^\9)((0,T)\times\rrd)\cap L^\9(0,T;L^2)$ will be proved for $s\in\left(\frac12,1\right)$ and $\beta'(r)>0$, $\ff\,r\in\rr$ and $\beta'\ge0$ if $D\equiv0$. In the special case of porous media equations with fractional Laplacian, that is, $D\equiv0$, $\beta(u)\equiv|u|^{m-1}u$, $m>(d-2s)_+/d$, the existence of a strong solution was proved in \cite{7}, \cite{7a} (see also \cite{6} for some earlier abstract results, which applies to this case as well). 

Like in the present work, the results obtained in \cite{7} are based on the Crandall \& Liggett generation theorem of nonlinear contraction semigroups in $L^1(\rrd)$. However, the approach used in \cite{7} cannot be adapted to cope with equation \eqref{e1.1}. In fact, the existence and uniqueness of a mild solution to \eqref{e1.1} reduces to prove the $m$-accretivity in $L^1(\rrd)$ of the operator $A_0$, that is, \mbox{$(I+\lbb A_0)\1$} must be nonexpansive in $L^1(\rrd)$ for all $\lbb>0$. If $D\equiv0$ and $\beta(u)=|u|^{m-1}u$, $m>(d-2s)_+/d$, this follows as shown in \cite{7} (see, e.g., Theorem 7.1) by regularity $u\in L^1(\rrd)\cap L^{m+1}(\rrd)$, $|u|^{m-1}u\in \dot H^s(\rrd)$ of solutions to the resolvent equation $u+\lbb(-\Delta)^s\beta(u)=f$ for $f\in L^1(\rrd)$. However, such a property might not be true in our case. For instance, if $s=1$, this happens if $|b'(r)r+b(r)|\le\alpha\beta'(r)$, $\ff r\in\rr,$ $b\ge0$, $\beta'>0$ on $\rr\setminus\{0\}$, and $D$ sufficiently regular (\cite[Theorem 2.2]{4}). To circumvent this situation, following \cite{2} (see Section \ref{s2}) we have constructed here an $m$-accretive restriction $A$ of $A_0$ and derive so via the Crandall \& Liggett theorem a semigroup of contractions $S(t)$ such that $u(t)=S(t)u_0$ is a mild solution to \eqref{e1.1}. In general, that is if $A\ne A_0$, this is not the unique mild solution to \eqref{e1.1}. However, as shown in Theorem \ref{t3.1} below, under Hypotheses (j) (resp. (j)$'$), (jj), (jjj) (see Section \ref{s3}), for initial conditions in $L^1\cap L^\9$ it is the unique bounded, distributional solution to \eqref{e1.1}. For initial conditions in $L^1$, the uniqueness of mild solutions to \eqref{e1.1} as happens for $s=1$ (\cite{4}) or for $D\equiv0$, $s\in(0,1)$, as shown in \cite{7}, in the case of the present paper remains open. One may suspect, however, that one has in this case as for $s=1$ (see \cite{8a}, \cite{8aa}) the existence   of an entropy, resp. kinetic, solution to \eqref{e1.1} for $u_0\in L^1\cap L^\9$. But this remains to be done. Let us mention that there is a huge literature on the well-posedness of equation \eqref{e1.1} for the case $s=1$, in particular when $D\equiv0$. We refer the reader e.g. to \cite{a2'}--\cite{4}, \cite{5}, \cite{8a}, \cite{8aa}, \cite{6}, \cite{a15'} and the references therein.

In Section \ref{s4}, we apply our results to the following  McKean--Vlasov SDE on $\rrd$
\begin{equation}
	\label{e1.11}\barr{l}
	dX_t=D(X_t)b(u(t,X_t))dt+\(\dd\frac{\beta(u(t,X_{t-}))}{u(t,X_{t-})}\)^{\frac1{2s}}dL_t,\vsp
	\call_{X_t}(dx):=\mathbb{P}\circ X^{-1}_t(dx)=u(t,x)dx,\ t\in[0,T],\earr
	\end{equation}where $L$ is a $d$-dimensional isotropic $2s$-stable process with L\'evy measure $dz/|z|^{d+2s}$ (see \eqref{e4.5} below).  We prove that provided $u(0,\cdot)$ is a probability density in $L^\9$, by our Theorem \ref{t2.3} and the superposition principle for non-local Kolmogorov operators (see \cite[Theorem 1.5]{7aa}, which is an extension of the local case in \cite{b17'} and \cite{13'})  it follows that \eqref{e1.11} has a weak solution (see Theorem \ref{t4.1} below). Furthermore, we prove that our Theorem \ref{t3.1} implies that we have weak uniqueness for \eqref{e1.11} among all solutions satisfying
$$\left((t,x)\mapsto\frac{d\call_{X_t}}{dx}\,(x)\right)\in L^\9((0,T)\times\rrd),$$(see Theorem \ref{t4.2}). As a consequence, their laws form a nonlinear Markov process in the sense of McKean \cite{9'}, thus realizing his vision formulated in that paper (see Remark 4.3). We stress that for the latter two results $\beta$ is allowed to be degenerate, if $D\equiv0$. We refer to Section \ref{s4} for details.

McKean--Vlasov SDEs for which $(L_t)$ in \eqref{e1.11} is replaced by a Wiener process $(W_t)$ have been studied very intensively following the two fundamental papers \cite{9'}, \cite{a23}. We refer to \cite{a11'}, \cite{a19'} and the monograph \cite{a9'} as well as the references therein. We stress that \eqref{e1.11} is of {\it Nemytskii type}, i.e. distribution density dependent, also called {\it singular} McKean--Vlasov SDEs, so there is no weak continuity in the measure dependence of the coefficients, as usually assumed in the literature. This (also in case of Wiener noise) is a technically more difficult situation. Therefore, the literature on weak existence and uniqueness for \eqref{e1.11} with L\'evy noise is much smaller. In fact, since the diffusion coefficient is allowed to 
depend (nonlinearly) on the distribution density, except for \cite{7aa}, where weak existence (but not uniqueness) is proved for \eqref{e1.11}, if $D\equiv0$ and $\beta(r):=|r|^{m-1}r,$ $m>(d-2\sigma)_+/d,$ we are not aware of any other paper adressing weak well-posedness in our case. If in \eqref{e1.11} the L\'evy process $(L_t)$ is replaced by a Wiener process $(W_t)$, we refer to \cite{a2'}, \cite{a2''}, \cite{a2'''}, \cite{2} for weak existence  and to \cite{3}--\cite{4} for weak uniqueness, as well as the references therein.

\bk\noindent{\bf Notation.} $L^p(\rrd)=L^p,\ p\in[1,\9]$ is the standard space of Lebesgue  $p$-integrable functions on $\rr^d$. We denote by $L^p_{\rm loc}$ the correspon\-ding local space and by $|\cdot|_p$ the norm of $L^p$. The inner product in $L^2$ is denoted by $(\cdot,\cdot)_2$. Denote by 
$H^\sigma(\rrd)=H^\sigma$, $0<\sigma<\9$, the standard Sobolev spaces on $\rrd$ in $L^2$ and by $H^{-\sigma}$ its dual space. By $C_b(\rr)$ denote the space of continuous and bounded functions on $\rr$ and by $C^1(\rr)$ the space of differentiable functions on~$\rr$. For any $T>0$ and a Banach space $\calx$, $C([0,T];\calx)$ is the space of $\calx$-valued continuous functions on $[0,T]$ and by $L^p(0,T;\calx)$ the space of $\calx$-valued $L^p$-Bochner integrable functions on $(0,T).$ We denote also by $C^\9_0(\calo)$, $\calo\subset\rrd$, the space of infinitely differentiable functions with compact support in $\calo$ and by $\cald'(\calo)$ its dual, that is, the space of Schwartz distributions on $\calo$. By $C^\9_0([0,\9)\times\rrd)$ we denote the space of infinitely differentiable functions on $[0,\9)\times\rrd$ with compact in $[0,\9)\times\rrd$. By $S'(\rrd)$ we denote the space of tempered distributions on $\rrd$.

\section{Existence of a mild solution}\label{s2}
\setcounter{equation}{0}

To begin with, let us construct the operator $A:D(A)\subset L^1\to L^1$ mentioned in \eqref{e1.4}.  To this purpose, we shall first prove  the following lemmas. 

\begin{lemma}\label{l2.1} Assume that $\frac12< s<1$. Let $\lbb_0>0$ be as defined in \eqref{a2.27''} below. Then, under Hypotheses {\rm(i)--(iv)} there is a family of operators $\{J_\lbb:L^1\to L^1;\lbb>0)\}$, which for all $\lbb\in(0,\lbb_0)$  satisfies
	\begin{eqnarray}
&	(I+\lbb A_0)(J_\lbb(f))=f,\ \ff\,f\in L^1,\label{e2.1}\\[1mm]
&	|J_\lbb(f_1)-J_\lbb(f_2)|_1\le|f_1-f_2|_1,\ \ff f_1,f_2\in L^1,\label{e2.2}\\[1mm]
&	J_{\lbb_2}(f)=J_{\lbb_1}\(\dd\frac{\lbb_1}{\lbb_2}\,f+\(1-\frac{\lbb_1}{\lbb_2}\)J_{\lbb_2}(f)\),\ \ff f\in L^1,\ \lbb_1,\lbb_2>0,\qquad\label{e2.3}
\\[1mm]
&\dd\int_\rrd J_\lbb(f)dx=\int_\rrd f\,dx,\ \ff f\in L^1,\label{e2.4}\\[1mm]
& J_\lbb(f)\ge0,\ \mbox{ a.e. on }\rrd,\mbox{ if }f\ge0,\mbox{ a.e. on }\rrd,\label{e2.5}\\[1mm]
& |J_\lbb(f)|_\9\le(1+||D|+({\rm div}\,D)^-|^{\frac12}_\9)|f|_\9,\ \ff\,f\in L^1\cap L^\9,	\label{e2.5a}\\[1mm]
&\beta(J_\lbb(f))\in H^{s}\cap L^1\cap L^\9,\ \ff f\in L^1\cap L^\9. \label{e2.6a}
\end{eqnarray}
\end{lemma}

\begin{remark}\label{r2.1'} \rm By \eqref{e2.5a}, \eqref{e2.6a} and for given $f\in L^1\cap L^\9$ changing $\beta$ as in the proof of Theorem \ref{t3.1} below, we may drop the assumption $\beta\in {\rm Lip}(\rr)$ in Hypothesis (i).\end{remark}

\n{\bf Proof of Lemma \ref{l2.1}.} We shall  first prove the existence of a solution $y=y_\lbb\in D(A_0)$ to the equation
\begin{equation}\label{e2.6} 
y+\lbb A_0(y)=f\mbox{ in }S',\end{equation}for $f\in L^1$. To this end, for $\vp\in(0,1]$ we consider the approximating equation
\begin{equation}\label{e2.7} 
y+\lbb (\vp  I-\Delta)^s(\beta_\vp(y))+\lbb\,{\rm div}(D_\vp b_\vp(y)y)=f\mbox{ in }S',\end{equation}where, for $r\in\rr$, $\beta_\vp(r):=\beta(r)+\vp r$ and
$$D_\vp:=\eta_\vp D,\ \eta_\vp\in C^1_0(\rrd),\ 0\le\eta_\vp\le1,\ |\nabla\eta_\vp|\le1,\ \eta_\vp(x)=1\mbox{ if }|x|<\frac1\vp.$$Clearly, we have
\begin{equation}
\label{e2.8a}
\barr{c}
|D_\vp|\in L^2\cap L^\9,\ |D_\vp|\le|D|,\ \dd\lim_{\vp\to\9}D_\vp(x)=D(x),\mbox{ a.e. }x\in \rrd,\vsp 
{\rm div}\,D_\vp\in L^2,\ ({\rm div}\,D_\vp)^-\le({\rm div}\,D)^-+\one_{\left[|x|>\frac1\vp\right]}|D|.\earr\end{equation}
As regards $b_\vp$, it is of the form 
$$b_\vp(r)\equiv \frac{(b*\vf_\vp)(r)}{1+\vp|r|},\ \ff\,r\in\rr,$$where 
 $\vf_\vp(r)=\frac1\vp\ \vf\(\frac r\vp\)$ is a standard mollifier.  We also set $b^*_\vp(r):=b_\vp(r)r,$ \mbox{$ r\in\rr.$} 
 
 Now, let us assume that $f\in L^2$ and consider the approximating equation
 \begin{equation}\label{2.11}
 	F_{\vp,\lbb}(y)=f\mbox{ in }S',\end{equation}where $F_{\vp,\lbb}:L^2\to S'$ is defined by
 $$F_{\vp,\lbb}(y):=y+\lbb(\vp I-\Delta)^s\beta_\vp(y)+\lbb\,{\rm div}(D_\vp b^*_\vp(y)),\ y\in L^2,$$where $(\vp I-\Delta)^s:S\to S$ is defined as usual by  Fourier transform and then it extends by duality to an operator $(\vp I-\Delta)^s:S'\to S'$ (which is consistent with \eqref{e1.2}). 
 
 We recall that the Bessel space of order $s\in\rr$ is defined as
 $$H^s:=\{u\in S';\ (1+|\xi|^2)^{\frac s2}\calf(u)\in L^2\}$$and the Riesz space as
 $$\dot H^s:=\{u\in S';\ \calf(u)\in L^1_{\rm loc}\mbox{ and }|\xi|^s\calf(u)\in L^2\}$$with respective norms
 $$|u|^2_{H^s} :=\int_\rrd(1+|\xi|^2)^s|\calf(u)|^2(\xi)d\xi
 =\dd\int_\rrd|(I-\Delta)^{\frac s2}u|^2d\xi,$$and 
 $$|u|^2_{\dot H^s} :=\int_\rrd|\xi|^{2s}|\calf(u)|^2(\xi)d\xi=\int_\rrd|(-\Delta)^{\frac s2}u|^2d\xi.$$
 
 $H^s$ is a Hilbert space for all $s\in\rr$, whereas $\dot H^s$ is only a Hilbert space if $s<\frac d2$ (see, e.g., \cite[Proposition 1.34]{0}). 
 
 Now, we shall show that \eqref{2.11} has a unique solution $y_\vp\in L^2$. To this end, we rewrite \eqref{2.11} as
 $$(\vp I-\Delta)^{-s}F_{\vp,\lbb}(y)=(\vp I-\Delta)^{-s}f\ (\in H^{2s}),$$i.e.,
 \begin{equation}
 	\label{2.12}
 	(\vp I-\Delta)^{-s}y+\lbb\beta_\vp(y)+\lbb(\vp I-\Delta)^{-s}{\rm div}(D_\vp b^*_\vp(y))=(\vp I-\Delta)^{-s}f.
 	\end{equation}
 Clearly, since $D_\vp b^*_\vp(y)\in L^2$, hence ${\rm div}(D_\vp b^*_\vp(y))\in H\1$, we have
 $$(\vp I-\Delta)^{-s}F_{\vp,\lbb}(y)\in L^2,\ \ff y\in L^2,$$because $s>\frac12$. Now, it is easy to see that \eqref{2.12} has a unique solution, $y_\vp\in L^2$, because, as the following chain of inequalities shows, $(\vp I-\Delta)^{-s}F_{\vp,\lbb}:L^2\to L^2$ is strictly monotone. By \eqref{2.12} we have, for $y_1,y_2\in L^2$,
\begin{equation}
	\label{2.13}
	\hspace*{-4mm}\barr{l}
 ((\vp I-\Delta)^{-s}(F_{\vp,\lbb}(y_2)-F_{\vp,\lbb}(y_1)),y_2-y_1)_2\vsp
 \quad=((\vp I-\Delta)^{-s}(y_2-y_1),y_2-y_1)_2
 +\lbb(\beta_\vp(y_2)-\beta_\vp(y_1),y_2-y_1)_2\vsp
 \quad-\lbb{\ }_{H^{-1}}\!\!\<{\rm div}(D_\vp(b^*_\vp(y_2)-b^*_\vp(y_1))),
 (\vp I-\Delta)^{-s}(y_2-y_1)\>_{H^1}\vsp
 \quad\ge|y_2-y_1|^2_{H^{-s}}+\lbb\vp|y_2-y_1|^2_2\vsp
 \qquad-\lbb c_1|D_\vp(b^*_\vp(y_2)-b^*_\vp(y_1))|_2|(\vp I-\Delta)^{\frac12-s}(y_2-y_1)|_2\vsp
 \quad\ge|y_2-y_1|^2_{H^{-s}}+\lbb\vp|y_2-y_1|^2_2-\lbb c|D|_\9{\rm Lip}(b^*_\vp)|y_2-y_1|_2|y_2-y_1|_{H^{1-2s}},
 \earr\hspace*{-6mm}\end{equation} 
where $c\in(0,\9)$ is independent of $\lbb,\vp,y_1,y_2$ and ${\rm Lip}(b^*_\vp)$ denotes the Lipschitz norm of $b^*_\vp$. Since $-s<1-2s<0$, by interpolation we have for $\theta:=\frac{2s-1}s$ that
$$|y_2-y_1|_{H^{1-2s}}\le|y_2-y_1|^{1-\theta}_2|y_2-y_1|^\theta_{H^{-s}}$$(see \cite[Proposition~1.52]{0}). So, by Young's inequality we find that the left hand side of \eqref{2.13} dominates
$$\lbb(\vp-\lbb c_\vp)|y_2-y_1|^2_2+\frac12|y_2-y_1|^2_{H^{-s}}$$for some $c_\vp\in(0,\9)$ independent of $\lbb,y_1$ and $y_2$. Hence, for some $\lbb_\vp\in(0,\9)$, we conclude that $(\vp I-\Delta)^{-s}F_{\vp,\lbb}$ is strictly monotone on $L^2$ for $\lbb\in(0,\lbb_\vp)$.

It follows from \eqref{2.12} that its solution $y_\vp$ belongs to $H^{2s-1}$, hence $b^*_\vp(y_\vp)\in H^{2s-1}$. Since $s>\frac12$ and $D\in C^1(\rrd;\rrd)$, by simple bootstrapping \eqref{2.12} implies
\begin{equation}
	\label{2.14}
	y_\vp\in H^1,\end{equation}hence $\beta_\vp(y_\vp)\in H^1$. Furthermore, for $f\in L^2$ and $\lbb\in(0,\lbb_\vp)$, $y_\vp$ is the unique solution of \eqref{e2.7}  in $L^2$. 
  
Assume now that $\lbb\in(0,\lbb_\vp)$ and $f\ge0$, a.e. on $\rrd$. 
Then, we have
\begin{equation}
\label{e2.13}
y_\vp\ge0,\ \mbox{ a.e. on }\rrd.
\end{equation}
Here is the argument. For $\delta>0$, consider the function
\begin{equation}
\label{e2.14}
\eta_\delta(r)=\left\{\barr{rll}
-1&\mbox{ for }&r\le-\delta,\vsp
\dd\frac r\delta&&r\in(-\delta,0),\vsp 
0&\mbox{ for }&r\ge0.\earr\right.
\end{equation}
If we multiply equation \eqref{e2.7}, where $y=y_\vp$, by $\eta_\delta(\beta_\vp(y_\vp))$ $(\in H^1)$ and integrate over $\rrd$, we get
\begin{equation}
\label{e2.15}
\barr{l}
\dd\int_\rrd y_\vp\eta_\delta(\beta_\vp(y_\vp))dx
+\lbb\dd\int_\rrd(\vp I-\Delta)^s(\beta_\vp(y_\vp))\eta_\delta(\beta_\vp(y_\vp))dx\vsp 
\qquad 
=\dd\int_\rrd f\eta_\delta(\beta_\vp(y_\vp))dx
+\lbb\dd\int_\rrd D_\vp b^*_\vp(y_\vp)\eta'_\delta(\beta_\vp(y_\vp))\cdot\nabla \beta_\vp(y_\vp) dx.
\earr\end{equation}By Lemma 5.2 in \cite{7} we have (Stroock-Varopoulos inequality) 
\begin{equation}
\label{e2.15a}
\int_\rrd(\vp I-\Delta)^su\Psi(u)dx\ge\int_\rrd|(\vp I-\Delta)^{\frac s2}\wt\Psi(u)|^2dx,\ u\in H^1(\rrd),
\end{equation}for any pair of functions $\Psi,\wt\Psi\in{\rm Lip}(\rr)$ such that $\Psi'(r)\equiv(\wt\Psi'(r))^2,\ r\in\rr.$ This yields
\begin{equation}\label{e2.15aa}
\dd\int_\rrd(\vp I-\Delta)^s\beta_\vp(y_\vp)\eta_\delta(\beta_\vp(y_\vp))dx
\ge\dd\int_\rrd|(\vp I-\Delta)^{\frac s2}\wt\Psi(\beta_\vp(y_\vp))|^2dx\ge0,\end{equation}where $\wt\Psi(r)=\int^r_0\sqrt{\eta'_\delta(s)}\,ds.$ 
Taking into account that $y_\vp\in H^1$ and that $|\beta_\vp(y_\vp)|\ge \vp|y_\vp|$, we have
\begin{equation}
\label{e2.16}
\barr{l}
\left|\dd\int_\rrd D_\vp b_\vp(y_\vp)y_\vp
\eta'_\delta(\beta_\vp(y_\vp))\nabla \beta_\vp(y_\vp) dx\right|\\
\qquad\le\dd\frac1\delta\,|b|_\9\int_{\wt E^\delta_\vp}|y_\vp(x)|\,|\nabla \beta_\vp(y_\vp)|\,|D_\vp(x)|dx\\
\qquad\le\dd\frac 1\vp\,|b|_\9\|D_\vp\|_{L^2}
\(\dd\int_{\wt E^\delta_\vp}|\nabla \beta_\vp(y_\vp)|^2dx\)^{\frac12}\to0\mbox{ as }\delta\to0.\earr
\end{equation}
Here  $\wt E^\delta_\vp=\{-\delta<\beta_\vp(y_\vp)\le0\}$ and we used that $\nabla\beta_\vp(y_\vp)=0$, a.e. on $\{\beta_\vp(y_\vp)=0\}$.  

Taking into account that ${\rm sign}\,\beta_\vp(r)\equiv{\rm sign}\,r$, by \eqref{e2.15}--\eqref{e2.16}  we get, for $\delta\to0$, that $y^-_\vp=0,$ a.e. on $\rrd$ and so \eqref{e2.13} holds.

If $y_\vp=y_\vp(\lbb,f)$ is the solution to \eqref{e2.7}, we have for $f_1,f_2\in L^1\cap L^2$
\begin{equation}
\label{e2.18}
\barr{l}
y_\vp(\lbb,f_1)-y_\vp(\lbb,f_2)
+\lbb(\vp I-\Delta)^s(\beta_\vp(y_\vp(\lbb,f_1))-\beta_\vp(y_\vp(\lbb,f_2)))\vsp\qquad+\lbb\ {\rm div}\ D_\vp(b^*_\vp(y_\vp(\lbb,f_1))-b^*_\vp(y_\vp(\lbb,f_2)))=f_1-f_2.\earr\end{equation}

Now, for $\delta>0$ consider the function
\begin{equation*}
\calx_\delta(r)=\left\{\barr{rll}
1&\mbox{ for }&r\ge\delta,\vsp
\dd\frac r\delta&&\mbox{ for }|r|<\delta,\vsp 
-1&\mbox{ for }&r<-\delta.\earr\right.
\end{equation*}
If we multiply \eqref{e2.18} by $\calx_\delta(\beta_\vp(y_\vp(\lbb,f_1))-\beta_\vp(y_\vp(\lbb,f_2)))$ $(\in H^1)$ and integrate over $\rrd$, we get
$$\barr{l}
\dd\int_\rrd(y_\vp(\lbb,f_1)-y_\vp(\lbb,f_2))\calx_\delta(\beta_\vp(y_\vp(\lbb,f_1))-\beta_\vp(y_\vp(\lbb,f_2)))dx\\ 
\quad\le\!\lbb\,\dd\frac1\delta\!\int_{E_\vp^\delta}\!|b^*_\vp(y_\vp(\lbb,f_1)){-}b^*_\vp(y_\vp(\lbb,f_2))|\,|D_\vp|
|\nabla(\beta_\vp(y_\vp(\lbb,f_1))
{-}\beta_\vp(y_\vp(\lbb,f_2)))|dx\vsp
\quad +|f_1-f_2|_1,\earr$$because,  by \eqref{e2.15a}, we have
$$\dd\int_{\rrd}(\vp I-\Delta)^s(\beta_\vp(y_\vp,f_1)-\beta_\vp(y_\vp,f_2))\calx_\delta(\beta_\vp(y_\vp,f_1)-\beta_\vp(y_\vp,f_2))dx\ge0.$$ Set $E^\delta_\vp=\{|\beta_\vp(y_\vp(\lbb,f_1))-\beta_\vp(y_\vp(\lbb,f_2))|\le\delta\}.$ 

Since $|\beta_\vp(r_1)-\beta_\vp(r_2)|\ge\vp|r_1-r_2|,\ D_\vp\in L^2(\rrd;\rrd)\mbox{ and }b^*_\vp\in {\rm Lip}(\rr),$ we~get that
$$\lim_{\delta\to0}\frac1\delta\!\int_{E^\delta_\vp}\!|b^*_\vp(y_\vp(\lbb,f_1)){-}b^*_\vp(y_\vp(\lbb,f_2))||D_\vp||\nabla(\beta_\vp(y_\vp(\lbb,f_1)){-}\beta_\vp(y_\vp(\lbb,f_2)))|dx{=}0,$$because $y_\vp(\lbb,f_i)\in H^1,\ i=1,2$, and 
$\nabla(\beta_\vp(y_\vp(\lbb,f_1))-\beta_\vp(y_\vp(\lbb,f_2)))=0$, a.e. on $\{\beta_\vp(y_\vp(\lbb,f_1))-\beta_\vp(y_\vp(\lbb,f_2))=0\}$. This yields
\begin{equation}
\label{e2.19}
|y_\vp(\lbb,f_1)-y_\vp(\lbb,f_2)|_1\le|f_1-f_2|_1,\ \ff\lbb\in(0,\lbb_\vp).
\end{equation}
In particular, it follows that
\begin{equation}
\label{e2.23a}
|y_\vp(\lbb,f)|_1\le|f|_1,\ \ff\,f\in L^1\cap L^2.
\end{equation}
Now let us remove the restriction $\lbb$ to be in $(0,\lbb_\vp)$. To this end define the operator $A_\vp:D_0(A_\vp)\to L^1$ by
$$\barr{rl}
A_\vp(y)\!\!\!&:=(\vp I-\Delta)^s(\beta_\vp(y))+{\rm div}(D_\vp b^*_\vp(y)),\vsp 
D_0(A_\vp)\!\!\!&:=\{y\in L^1;\ (\vp I-\Delta)^s\beta_\vp(y)+{\rm div}(D_\vp b^*_\vp(y))\in L^1\},\earr$$and for $\lbb\in(0,\lbb_\vp)$
$$J^\vp_\lbb(f):=y_\vp(\lbb,f),\ f\in L^1\cap L^2.$$
Then $J^\vp_\lbb(L^1\cap L^2)\subset D_0(A_\vp)\cap H^1$ and by \eqref{e2.19} it extends by continuity to an operator $J^\vp_\lbb:L^1\to L^1$. We note that the operator $(A_\vp,D_0(A_\vp))$ is closed as an operator on $L^1$. Hence \eqref{e2.19} implies that
\begin{equation}\label{2.22'}
	J^\vp_\lbb(L^1)\subset D_0(A_\vp)\end{equation}and that $J^\vp_\lbb(f)$ solves \eqref{e2.7} for all $f\in L^1$. 

Now define
\begin{equation}\label{2.22''}
	D(A_\vp):=J^\vp_\lbb(L^1)\end{equation}
and restrict $A_\vp$ to $D(A_\vp)$. It is easy to see that $D(A_\vp)$ is independent of $\lbb\in(0,\lbb_\vp)$.

Now let $0<\lbb_1<\lbb_\vp$. Then, for $\lbb\ge\lbb_\vp$, the equation
\begin{equation}
	\label{e2.23aaaa}
	y+\lbb A_\vp(y)=f  \ (\in L^1),\ y\in D(A_\vp),
\end{equation}can be rewritten as

$$y+\lbb_1 A_\vp(y)=\(1-\frac{\lbb_1}\lbb\)y+\frac{\lbb_1}\lbb\,f.$$Equivalently,
\begin{equation}
	\label{e2.23v}
	y=J^\vp_{\lbb_1}\(\(1-\frac{\lbb_1}\lbb\)y+\frac{\lbb_1}\lbb\,f\). 
\end{equation}
Taking into account that, by \eqref{e2.19}, 
$|J^\vp_{\lbb_1}(f_1)-J^\vp_{\lbb_1}(f_2)|_1\le|f_1-f_2|_1,$ it follows that \eqref{e2.23v} has a unique solution $y_\vp\in D(A_\vp)$. Let  $J^\vp_\lbb(f):=y_\vp$, $\lbb>0,\ f\in L^1$, denote this solution to \eqref{e2.23aaaa}. By \eqref{e2.23v} we see that \eqref{e2.19}, \eqref{e2.23a} extend to all $\lbb>0,$  $f\in L^1.$ 

\mk\n{\bf Claim.} Let $f\in L^1\cap L^2$. Then
\begin{equation}\label{b2.27'}
J^\vp_\lbb(f)\in H^1  \mbox{ for all } \lbb>0.\end{equation}

\n{\bf Proof.} Fix $\lbb_1\in[\lbb_\vp/2,\lbb_\vp)$ and set $\lbb:=2\lbb_1.$ Define $T:L^1\cap L^2\to L^1\cap H^1$~by
$$T(y):=J^\vp_{\lbb_1}\(\frac12\ y+\frac12\ f\),\ y\in L^1\cap L^2.$$
Then, as just proved, for any $f_0\in L^1\cap L^2 $ fixed
\begin{equation}\label{b2.27''}
	\lim_{n\to\9}T^n(f_0)=J^\vp_\lbb(f)\mbox{ in }L^1.\end{equation}It suffices to prove
\begin{equation}\label{b2.27'''}
	J^\vp_\lbb(f)\in L^2  ,\end{equation}
because then $J^\vp_\lbb(f)=J^\vp_{\lbb_1}(g)$ with $g:=\frac12\, J^\vp_{\lbb}(f)+\frac12\ f\in L^1\cap L^2,$ and so the claim follows by \eqref{2.14} which holds with $y_\vp:=J^\vp_{\lbb_1}(g)$, because $\lbb_1\in(0,\lbb_\vp).$

To prove \eqref{b2.27'''} we note that  we have, for $n\in\nn$,
$$(I+\lbb_1 A_\vp)T^n(f_0)=\frac12\ T^{n-1}(f_0)+\frac12\ f$$with $T^n(f_0)\in H^1.$ Hence, applying ${}_{H^{-1}}\!\!\<\cdot,T^n(f_0)\>_{H^1}$ to this equation, we find
\begin{equation}\label{b2.27iv}
\!\!\barr{l}
|T^{(n)}f_0|^2_2+\lbb_1\, 
{}_{H^{-1}}\!\!\<(\vp I-\Delta)^s\beta_\vp(T^n(f_0)),T^n(f_0)\>_{H^1}\vsp 
 =\lbb_1\dd\int_\rrd(D_\vp b^*_\vp(T^n(f_0)))\cdot\nabla(T^n(f_0))d\xi
{+}\(\dd\frac12T^{n-1}(f_0)+\frac12f,T^n(f_0)\!\)_2\!\!.\earr\!\!\!\end{equation}
Setting
\begin{equation}\label{b2.27v}
\psi_\vp(r):=\int^r_0 b^*_\vp(\tau)d\tau,\ r\in\rr,\end{equation}
by Hypotheses (iii), (iv) we have
$$0\le\psi_\vp(r)\le|b^*_\vp|_\9r,\ r\in\rr,$$and hence the left hand side of \eqref{b2.27iv} is equal to
$$-\lbb_1({\rm div}\ D_\vp,\psi_\vp(T^n(f_0)))_2+\(\frac12\,(T^{n-1}(f_0)+f),T^n(f_0)\)_2,$$where we recall that ${\rm div}\,D_\vp\in L^2$ by \eqref{e2.8a}. 
By \eqref{e2.15aa} we thus obtain
$$\barr{ll}
|T^{(n)}(f_0)|^2_2\le\!\!&\lbb_1|b^*_\vp|_\9|({\rm div}\,D_\vp)^-|_2|T^n(f_0)|_2
+\dd\frac12|T^n(f_0)|^2_2\vsp &+\dd\frac14(|T^{n-1}(f_0)|^2_2+|f|^2_2),
\earr$$therefore,
$$|T^{(n)}(f_0)|^2_2\le C_\vp+\frac12\,|T^{n-1}(f_0)|^2_2,$$where
$$C_\vp:=\frac{\lbb_1}{1-\lbb_1}\,|b^*_\vp|^2_\9|({\rm div}\,D_\vp)^-|^2_2+\frac12\,|f|^2_2.$$Iterating, we find
$$|T^{(n)}(f_0)|^2_2\le \sum^{n-2}_{k=1} 2^{-k}C_\vp+
\frac1{2^{n-1}}\,|T(f_0)|^2_2,\ n\in\nn.$$
Hence, by Fatou's lemma and \eqref{b2.27''},
$$|J^\vp_\lbb(f)|^2_2\le\liminf_{n\to\9}|T^n(f_0)|^2_2\le C_\vp<\9,$$and \eqref{b2.27'} holds for $\lbb=2\lbb_1.$ Proceeding this way, we get \eqref{b2.27'} for all $\lbb>0.$~\hfill$\Box$\bk

Set
\begin{equation}\label{a2.27''}
\lbb_0:=\(\left|({\rm div}\ D)^-+|D|\right|_\9\)+
\(\left|({\rm div}\ D)^-+|D|\right|_\9^{\frac12}|b|_\9\)\1.\end{equation}
Then, for $f\in L^1\cap L^\9$ and $y_\vp:=J^\vp_\lbb(f),$
$\lbb>0$, we have\begin{equation}
\label{e2.23a2}
|y_\vp|_\9\le(1+\left||D|+({\rm div}\,D)^-\right|^{\frac12}_\9)|f|_\9,\ \ff\lbb\in(0,\lbb_0). 
\end{equation}
Indeed, if we set $M_\vp=|({\rm div}\,D_\vp)^-|^{\frac12}_\9|f|_\9,$  we get by \eqref{e2.7} that
$$\barr{c}
(y_\vp-|f|_\9-M_\vp) 
+\lbb(\vp I-\Delta)^s(\beta_\vp(y_\vp)-\beta_\vp(|f|_\9+M_\vp))
+\vp^s(\beta_\vp(|f|_\9+M_\vp))\vsp\qquad
+\lbb\ {\rm div}(D_\vp(b^*_\vp(y_\vp)-b^*_\vp(|f|_\9+M_\vp)))\le0.\earr$$Here we used that
$$1\in\{u\in S';\ (\vp+|\vp|^2)^s\calf(u)\in S'\}$$and that $(\vp I-\Delta)^s1=\vp^s,$ since $\calf(1)=\delta_0$ (= Dirac measure in $0\in\rrd$).
Then, taking the scalar product in $L^2$ with $\calx_\delta((\beta_\vp(y_\vp)-\beta_\vp(|f|_\9+M_\vp)^+)$, letting $\delta\to0$ and using \eqref{e2.15a},  we get by \eqref{e2.8a} 
$$y_\vp\le(1+||D|+({\rm div}\,D)^{-}|^{\frac12}_\9)|f|_\9,\mbox{ a.e. in }\rrd,$$and, similarly, for $-y_\vp$ which yields \eqref{e2.23a2} for $\lbb\in(0,\lbb_0)$.  

In particular, it follows that
\begin{equation}
\label{e2.23aa}
|J^\vp_\lbb(f)|_1+|J^\vp_\lbb(f)|_\9\le c_1,\ \ff\vp,\lbb>0,  
\end{equation}where $c_1=c_1(|f|_1,|f|_\9)$ is independent of $\vp$ and $\lbb$.

	Now, fix $\lbb\in(0,\lbb_0)$ and $f\in L^1\cap L^\9$. For $\vp\in(0,1]$ set
	$$y_\vp:=J^\vp_\lbb(f).$$
Then, since $\beta_\vp(y_\vp)\in H^1,$ by \eqref{e2.7} we get
\begin{equation}\label{e2.27'}
	\barr{l}
	(y_\vp,\beta_\vp(y_\vp))_2+\lbb{\ }_{H^{-1}}\!\<(\vp I-\Delta)^s\beta_\vp(y_\vp),\beta_\vp(y_\vp)\>_{H^1}\vsp 
	\qquad=-\lbb
	({\rm div}(D_\vp b^*_\vp(y_\vp)),\beta_\vp(y_\vp))_2
		+(f,\beta_\vp(y_\vp))_2\vsp 
		\qquad=\lbb\dd\int_\rrd(D_\vp b^*_\vp(y_\vp))\cdot\nabla\beta_\vp(y_\vp)dx+(f,\beta_\vp(y_\vp))_2.\earr\end{equation}

Setting 
\begin{equation}
\label{e2.36'}
\wt\psi_\vp(r):=\int^r_0 b^*_\vp(\tau)\beta'_\vp(\tau)d\tau,\ r\in\rr,
\end{equation}by Hypotheses (iii), (iv)  we have 
$$0\le\wt\psi_\vp(r)\le\frac12\ |b|_\9(|\beta'|_\9+1) r^2,\ \ff r\in\rr,$$and hence, since $y_\vp\in H^1$,  the left hand side of \eqref{e2.27'} is equal to
$$-\lbb\int_\rrd{\rm div}\,D_\vp\wt\psi_\vp(y_\vp)dx+(f,\beta_\vp(y_\vp))_2,$$which, because $(y_\vp,\beta_\vp(y_\vp))_2\ge0$ and $H^1\subset H^s$, by \eqref{e2.8a} and Hypothesis (iv) implies that
$$\lbb|(\vp I-\Delta)^{\frac s2}\beta_\vp(y_\vp)|^2_2
\le\frac12\,(\lbb|b|_\9(|\beta'|_\9+1)
|({\rm div}\,D)^-+|D||_\9)|y_\vp|^2_2
+\frac12|\beta_\vp(y_\vp)|^2_2+\frac12\,|f|^2_2.$$Since $|\beta_\vp(r)|\le({\rm Lip}(\beta)+1)|r|,\ r\in\rr,$ by \eqref{e2.23a}, \eqref{e2.23a2}   we obtain 
\begin{equation}\label{e2.27''}
	\sup_{\vp\in(0,1]}|(\vp I-\Delta)^{\frac s2}\beta_\vp(y_\vp)|^2_2\le C<\9,
	\end{equation}
where
$$C:=\lbb(|b|_\9+1)(|({\rm div}\,D)^-+|D|_\9|+2)^2({\rm Lip}(\beta)+1)^2|f|_\9|f|_1.$$
Since obviously for all $u\in H^s\ (\subset \dot H^s)$, $\vp\in(0,1]$,
\begin{equation}\label{b2.33'}
\!\!	|(-\Delta)^{\frac s2}u|^2_{2}\le|(\vp I-\Delta)^{\frac s2}u|^2_{2}
\le|(-\Delta)^{\frac s2}u|^2_{2}{+}\vp^s|u|^2_{2}
\!\le2|(\vp I{-}\Delta)^{\frac s2}u|^2_{2},\!\!\end{equation}and since $\beta_\vp(y_\vp)\in H^1\subset H^s$ with $|\beta_\vp(r)|\le(1+{\rm Lip}(\beta))|r|,$ we conclude from \eqref{e2.23aa}  and \eqref{e2.27''} that (along a subsequence) as $\vp\to0$
$$\barr{rcl}
\beta_\vp(y_\vp)&\to& z\mbox{ weakly in }\dot H^s,\vsp 
(\vp I-\Delta)^s\beta_\vp(y_\vp)&\to&(-\Delta)^sz\mbox{ in }S',\vsp 
y_\vp&\to&y\mbox{ weakly in $L^2$ and weakly$^*$ in $L^\9$},\earr$$
where the second statement follows, because 
$$(\vp I-\Delta)^s\vf\to(-\Delta)^s\vf\mbox{ in }L^2\mbox{ for all }\vf\in S.$$By \cite[Theorem 1.69]{0}, it follows that as $\vp\to0$
$$\beta_\vp(y_\vp)\to z\mbox{ in }L^2_{\rm loc}(\rrd),$$so (selecting another subsequence, if necesary)
$$\beta(y_\vp)\to z,\mbox{ a.e. }$$
Since $\beta\1$ (the inverse function of $\beta$) is continuous, it follows that as $\vp\to0$
$$y_\vp\to\beta\1(z)=y,\mbox{ a.e.},$$so
\begin{equation}\label{2.32'}
	z=\beta(y)\end{equation}and
$$b^*_\vp(y_\vp)\to b^*(y)\mbox{ weakly in }L^2.$$Recalling that $y_\vp$ solves \eqref{e2.7}, we can let $\vp\to0$ in \eqref{e2.7} to find that 
\begin{equation}\label{2.27'''}
	y+\lbb(-\Delta)^s\beta(y)+\lbb\,{\rm div}(Db^*(y))=f\mbox{ in }S'.\end{equation}Since $\beta\in {\rm Lip}(\rr)$, the operator $(A_0,D(A_0))$ defined in \eqref{e1.5}  is obviously closed as an operator on $L^1$. Again defining for $y$ as in \eqref{2.27'''}
$$J_\lbb(f):=y\in D(A_0),$$  it follows by \eqref{e2.19} and Fatou's lemma that for $f_1,f_2\in L^1\cap L^\9$
\begin{equation}\label{2.27iv}
	|J_\lbb(f_1)-J_\lbb(f_2)|_1\le|f_1-f_2|_1.\end{equation}
Hence $J_\lbb$ extends continuously to all of $L^1$, still satisfying \eqref{2.27iv} for all $f_1,f_2\in L^1$. Then it follows by the closedness of $(A_0,D(A_0))$ on $L^1$ that $J_\lbb(f)\in D(A_0)$ and that it solves \eqref{2.27'''} for all $f\in L^1$. 

Hence, Lemma \ref{l2.1} is proved except for \eqref{e2.3} and \eqref{e2.4}. However, \eqref{e2.3} is obvious, since by \eqref{e2.1} it is equivalent to
$$(I+\lbb_1A_0)J_{\lbb_2}(f)=\frac{\lbb_1}{\lbb_2}\,f+\(1-\frac{\lbb_1}{\lbb_2}\)J_{\lbb_2}f,$$which in turn is equivalent to
$$(I+\lbb_2A_0)J_{\lbb_2}(f)=f.$$
Now let us prove \eqref{e2.4}. We may assume that $f\in L^1\cap L^\9$ and set $y:=J_\lbb(f)$. Let $\calx_n\in C^\9_0(\rrd),\ \calx_n\uparrow1$, as $n\to\9$, with $\sup\limits_n|\nabla\calx_n|_\9<\9$. Define
$$\vf_n:=(I+(-\Delta)^s)\1\calx_n=g^s_1*\calx_n,\ n\in\nn,$$where $g^s_1$ is as in the Appendix. Then, clearly,
\begin{equation}
	\label{b2.36}
	\barr{l}
	\vf_n\uparrow1,\mbox{ as }n\to\9,\vsp 
	\sup\limits_n(|\vf_n|_\9+|\nabla\vf_n|_\9)<\9,\vsp 
	\vf_n\in L^1\cap H^{2s},\ n\in\nn,\earr
	\end{equation}
where the last statement follows from \eqref{b2.33'}. Furthermore,
$$(-\Delta)^s\vf_n=\calx_n-(I+(-\Delta)^s)\1\calx_n\in L^1\cap L^\9,$$and, as $n\to\9$,
$$(-\Delta)^s\vf_n\to0\ dx-\mbox{a.e.},$$hence, because $\beta(y)\in L^1\cap L^\9,$
\begin{equation}\label{b2.37}
	\lim_{n\to\9}\int_\rrd(-\Delta)^s\vf_n\,\beta(y)dx=0.\end{equation}
Consequently, since $\beta(y)\in H^s$, $y\in D(A)$ with $A_0y\in L^1$,
$$\barr{l}
\dd\int_\rrd A_0y\,dx\vsp 
=\dd\lim_{n\to\9}\int_\rrd\vf_n A_0y\,dx\vsp
=-\dd\int_\rrd\beta(y)dx+\lim_{n\to\9}
\,{}_{H^{2s}}\!\!\<\vf_n,(I+(-\Delta)^s)\beta(y)+{\rm div}(Db^*(y))\>_{H^{-2s}}\vsp 
=-\dd\int_\rrd\!\beta(y)dx{+}\lim_{n\to\9}\, {}_{H^s}\!\!\<\vf_n,(I{+}(-\Delta)^s)\beta(y)\>_{H^{-s}}
\!+\!\dd\lim_{n\to\9}\int_\rrd\!\!\nabla\vf_n\cdot Db^*(y)dx,
\earr$$which by \eqref{b2.36} and \eqref{b2.37} is equal to zero. Hence, integrating the equation
$$y+\lbb A_0y=f$$over $\rrd$, \eqref{e2.4} follows, which concludes the proof of Lemma \ref{l2.1}.~\hfill$\Box$\bk

Now define
\begin{equation}\label{2.35}
	\barr{rl}
		D(A)\!\!\!&:=J_\lbb(L^1)\ (\subset D(A_0)),\vsp 
		A(y)\!\!\!&:=A_0(y),\ y\in D(A).\earr\end{equation}
	 Again it is easy to see that $J_\lbb(L^1)$ is independent of $\lbb\in(0,\lbb_0)$ and that
$$J_\lbb=(I+\lbb A)\1,\ \lbb\in(0,\lbb_0).$$
Therefore, we have
\begin{lemma}\label{l2.2} Under Hypotheses  {\rm(i)--(iv)}, the operator $A$ defined by \eqref{2.35} is $m$-accretive in $L^1$ and $(I+\lbb A)\1=J_\lbb$, $\lbb\in(0,\lbb_0)$. Moreover, if $\beta\in C^\9(\rr)$, then $\ov{D(A)}=L^1$.\end{lemma}
Here, $\ov{D(A)}$ is the closure of $D(A)$ in $L^1$.  

We note that, by \eqref{e1.3}, if $\beta\in C^\9(\rr)$, it follows that
$$A_0(\vf)=(-\Delta)^s\beta(\vf)+{\rm div}(Db(\vf)\vf)\in L^1,\ \ff\vf\in C^\9_0(\rrd),$$and so $\ov{D(A)}=L^1$, as claimed. 

Then, by  the Crandall \& Liggett theorem (see, e.g., \cite{1}, p.~131), we have that the Cauchy problem \eqref{e1.4}, that is,
$$\barr{l}
\dd\frac{du}{dt}+A(u)=0,\ \ t\ge0,\vsp
u(0)=u_0,\earr$$has, for each $u_0\in\ov{D(A)}$, a unique mild solution $u=u(t,u_0)\in C([0,\9);L^1)$ and $S(t)u_0=u(t,u_0)$ is a $C_0$-semigroup of contractions on $L^1$, that is,
$$\barr{c}
|S(t)u_0-S(t)\bar u_0|_1\le|u_0-\bar u_0|_1,\ \ff u_0,\bar u_0\in\ov{D(A)},\vsp 
S(t+\tau)u_0=S(t,S(\tau)u_0),\ \ff t,\tau>0;\ u_0\in\ov{D(A)},\vsp 
\dd\lim_{t\to0}S(t)u_0=u_0\mbox{\ \ in }L^1(\rrd).\earr$$  
Moreover, by \eqref{e2.5} and the exponential formula
$$S(t)u_0=\lim_{n\to\9}\(I+\frac tn\,A\)^{-n}u_0,\ \ff\,t\ge0,$$it follows that $S(t)u_0\in L^\9((0,T)\times\rrd)$, $T>0$, if $u_0\in L^1\cap L^\9.$ 

Let us show now that $u=S(t)u_0$ is a Schwartz distributional solution, that is, \eqref{e1.10} holds.

By \eqref{e1.6}--\eqref{e1.9}, we have
$$\barr{l}
\dd\int^\9_0dt\(\int_\rrd\vf(t,x)
(u_h(t,x)-u_h(t-h,x)\)dx\vsp
\qquad+\dd\int_\rrd
(\vf(t,x)(-\Delta)^s\beta(u_h(t,x))
-\nabla_x\vf(t,x)\cdot D(x)b^*(u_h((x)))dx)=0,\\\hfill \ff\vf\in C^\9_0([0,\9)\times\rrd).\earr$$This yields
$$\barr{l}
\dd\frac1h\int^\9_0dt
\(\int_\rrd u_h(t,x)(\vf(t+h,x)-\vf(t,x))\)dx\vsp
 \qquad+\dd\int_\rrd(\beta(u_h(t,x))(-\Delta)^s\vf(t,x)
 -\nabla_x\vf(t,x)\cdot D(x)b^*(u_h(t,x))dx)\vsp 
\qquad +\dd\frac1h\int^h_0dt\int_\rrd u_0(x)\vf(t,x)dx=0,\ \ff\vf\in C^\9_0([0,\9)\times\rrd).
\earr$$Taking into account that, by \eqref{e1.6} and (i)--(iii), $\beta(u_h)\to\beta(u)$, $b^*(u_h)\to b^*(u)$ in $C([0,T];L^1)$ as $h\to0$ for each $t>0$, we get  that \eqref{e1.10} holds. 
 
This together with Remark \ref{r2.1'} implies the following existence result for equation~\eqref{e1.1}.

 \begin{theorem}\label{t2.3} Assume that Hypotheses {\rm(i)--(iv)} hold. Then, there is a $C_0$-semigroup of contractions $S(t):L^1\to L^1$, $t\ge0$, such that for each \mbox{$u_0\in \ov{D(A)}$,} which is $L^1$ if $\beta\in C^\9(\rr)$, $u(t,u_0)=S(t)u_0$ is a mild  solution to \eqref{e1.1}. Moreover, if $u_0\ge0$, a.e. in~$\rrd$,
 	\begin{equation}
 	\label{e2.28}
 	u(t,u_0)\ge0,\mbox{\ \ a.e. in }\rrd,\ \ff\,t\ge0,
 	\end{equation}and
 	\begin{equation}
 	\label{e2.29}
 	\int_\rrd u(t,u_0)(x)dx=\int_\rrd u_0(x)dx,\ \ \ff\,t\ge0.
 	\end{equation}
 
 Moreover, $u$ is a distributional solution to \eqref{e1.1} on $[0,\9)\times\rrd.$ Finally, if $u_0\in L^1\cap L^\9$,  then all above assertions remain true, if we drop the assumption $\beta\in{\rm Lip}(\rr)$ from Hypothesis {\rm(i)}, and additionally we have that \mbox{$u\in L^\9((0,T)\times\rrd)$,} $T>0$.\end{theorem}

      \begin{remark}\label{r2.4}\rm It should be emphasized that, in general, the mild solution $u$ given by Theorem \ref{t2.3} is not unique because the operator $A$ constructed in Lemma \ref{l2.2} is dependent of the approximating operator $A_\vp y\equiv(\vp I+(-\Delta)^s)\beta_\vp(y)+{\rm div}(D_\vp b_\vp(y)y)$ and so $u=S(t)u_0$ may be viewed as a {\it viscosity-mild} solution to \eqref{e1.1}. However, as seen in the next section, this mild solution -- which is also a distributional solution to \eqref{e1.1} -- is, under appropriate assumptions on $\beta$, $D$ and $b$, unique in the class of solutions $u\!\in\! L^\9((0,T){\times}\rrd)$, $T>0$.\end{remark}
      
\section{The uniqueness of   distributional  solutions}\label{s3}
\setcounter{equation}{0}

In this section, we  shall prove   the uniqueness of distributional solutions to \eqref{e1.1}, where $s\in\(\frac12,1\)$, under the following Hypotheses: 
\bit\item[\rm(j)] $\beta\in C^1(\rr),\ \beta'(r)>0,\ \ff\,r\in\rr,\ \beta(0)=0.$
\item[\rm(jj)] $D\in L^\9(\rrd;\rrd).$
\item[\rm(jjj)] $b\in C^1(\rr).$
\eit 

\begin{theorem}\label{t3.1} Let $d\ge1,$   $s\in\left[\frac12,1\right]$, $T>0$,  and let $y_1,y_2\in L^\9((0,T){\times}\rrd)$ be two distributional solutions to \eqref{e1.1} on $(0,T)\times\rr^d$ $($in the sense of \eqref{e1.10}$)$   such that  $y_1-y_2\in L^1((0,T)\times\rrd)\cap L^\9(0,T;L^2)$ and 
\begin{equation}
\label{e3.1}\lim_{t\to0}\ {\rm ess} \sup_{\hspace*{-4mm}s\in(0,t)}|(y_1(s)-y_2(s),\vf)_2|=0,\ \ff\vf\in C^\9_0(\rrd).
\end{equation}	Then $y_1\equiv y_2$. If $D\equiv0$, then Hypothesis {\rm(j)} can be relaxed to
\begin{itemize}
	\item[\rm(j)$'$] $\beta\in C^1(\rr),\ \beta'(r)\ge0,\ \ff\,r\in\rr,\ \beta(0)=0.$
\end{itemize}\end{theorem}

\n{\bf Proof.} (The idea of proof is borrowed from Theorem 3.2 in \cite{4}, but has to be adapted substantially.) Replacing, if necessary, the functions $\beta$ and $b$ by
$$\beta_N(r)=\left\{\barr{ll}
\beta(r)&\mbox{ if }|r|\le N,\vsp
\beta'(N)(r-N)+\beta(N)&\mbox{ if }r>N,\vsp
\beta'(-N)(r+N)+\beta(-N)&\mbox{ if }r<-N,\earr\right.$$ and 
$$ b_N(r)=\left\{\barr{ll}
 b(r)&\mbox{ if }|r|\le N,\vsp
 b'(N)(r-N)+b(N)&\mbox{ if }r>N,\vsp
 b'(-N)(r+N)+b(-N)&\mbox{ if }r<-N,\earr\right.$$where $N\ge\max\{|y_1|_\9,|y_2|_\9\}$, by (j) we may assume that 
\begin{equation}\label{e3.2}
\beta',b'\in C_b(\rr),\ \beta'>\alpha_2\in(0,\9)\end{equation}and, therefore,   we have

\begin{eqnarray}
\alpha_1|\beta(r)-\beta(\bar r)|&\ge& 
|b^*(r)-b^*(\bar r)|,\ \ \ff\,r,\bar r\in\rr,\label{e3.3}\\[1mm]
|\beta(r)-\beta(\bar r)|&\ge&\alpha_2
|r-\bar r|,\ \ \ff\,r,\bar r\in\rr,\label{e3.3a}\end{eqnarray}
where $b^*(r)=b(r)r$ and $\alpha_2>0$. We set
\begin{equation}\label{e3.4}
\barr{c}
\Phi_\vp(y)=(\vp I+(-\Delta)^s)\1y,\ \ff\,y\in L^2,\vsp
z=y_1-y_2,\ w=\beta(y_1)-\beta(y_2),\ b^*(y_i)\equiv b(y_i)y_i,\ i=1,2.\earr\end{equation} 
It is well known that $\Phi_\vp:L^p\to L^p$, $\ff p\in[1,\9]$ and
\begin{equation}\label{e3.5}
\vp|\Phi_\vp(y)|_p\le |y|_p,\ \ \ff y\in L^p,\ \vp>0.
\end{equation}
Moreover, $\Phi_\vp(y)\in C_b(\rr^d)$ if $y\in L^1\cap L^\9$. We have 
$$z_t+(-\Delta)^s w+{\rm div}\,  D(b^*(y_1)-b^*(y_2))=0\mbox{ in }\cald'((0,T)\times\rrd).$$
We set
$$z_\vp=z*\theta_\vp,\ w_\vp=w*\theta_\vp,\ \zeta_\vp=(D(b^*(y_1)-b^*(y_2)))*\theta_\vp,$$where $\theta\in C^\9_0(\rr^d),$ $\theta_\vp(x)\equiv\vp^{-d}\theta\(\frac x\vp\)$ is a standard mollifier. We note that $z_\vp,w_\vp,\zeta_\vp,(-\Delta)^sw_\vp,{\rm div}\,\zeta_\vp\in L^2(0,T;L^2)$ and  we have
\begin{equation}
\label{e3.6}
(z_\vp)_t+(-\Delta)^s w_\vp+{\rm div}\,\zeta_\vp=0\mbox{ in }\cald'(0,T;L^2) .\end{equation}
This yields  $\Phi_\vp(z_\vp),\Phi_\vp(w_\vp),{\rm div}\,\Phi_\vp(\zeta_\vp)\in L^2(0,T;L^2)$ and 
\begin{equation}\label{e3.7}\barr{ll}
(\Phi_\vp(z_\vp))_t=-(-\Delta)^s\Phi_\vp(w_\vp)-{\rm div}\Phi_\vp(\zeta_\vp)=0\mbox{ in }\cald'(0,T;L^2).\earr\end{equation}
By \eqref{e3.6}, \eqref{e3.7} it follows that $(z_\vp)_t=\frac d{dt}\,z_\vp$,   $(\Phi_\vp(z))_t=\frac d{dt}\,\Phi_\vp(z_\vp)\in L^2(0,T;L^2)$, where $\frac d{dt}$ is taken in the sense of $L^2$-valued vectorial distributions on $(0,T)$. This implies that $z_\vp, \Phi_\vp(z_\vp)\in H^1(0,T;L^2)$ and both $[0,T]\ni t\mapsto z_\vp(t)\in L^2$ and $[0,T]\ni t\to\Phi_\vp(z_\vp(t))\in L^2$ are absolutely continuous. As a matter of fact, it follows by \eqref{e3.5} and \eqref{e3.7}  that
\begin{equation}\label{e3.8}
\Phi_\vp(z_\vp),\Phi_\vp(w_\vp)\in L^2(0,T;C_b(\rr^d)\cap L^2).\end{equation}
We set $h_\vp(t)=(\Phi_\vp(z_\vp(t)),z_\vp(t))_2$ and get, therefore,
\begin{eqnarray}
h'_\vp(t)\!&\!\!=\!\!&\!2(z_\vp(t),(\Phi_\vp(z_\vp(t)))_t)_2\label{e3.9}\\
\!&\!\!=\!\!&\!2(\vp\Phi_\vp(w_\vp(t)){-}w_\vp(t){-}{\rm div}\Phi_\vp(\zeta_\vp(t)),z_\vp(t))_2\nonumber\\
\!&\!\!=\!\!&\!2\vp(\Phi_\vp(z_\vp(t)),w_\vp(t))_2{+}2(\nabla\Phi_\vp(z_\vp(t)),\zeta_\vp(t))_2\nonumber\\
&&-2(z_\vp(t),w_\vp(t))_2,\mbox{ a.e. }t\in(0,T),\nonumber
\end{eqnarray}where, $(\cdot,\cdot)_2$ is the scalar product in $L^2$. 
By \eqref{e3.7}--\eqref{e3.9} it follows that \mbox{$t\to h_\vp(t)$} has an absolutely continuous $dt$-version on $[0,T]$ which we shall consider from now on. 
Since, by \eqref{e3.4}, we have for $\alpha_3:=|\beta'|^{-1}_\9$
\begin{equation}\label{e3.10}
(z_\vp(t),w_\vp(t))_2\ge\alpha_3|\,|w(t)|*\theta_\vp|^2_2+\gamma_\vp(t),
\end{equation}where
\begin{equation}
\label{e3.11}
\gamma_\vp(t):=(z_\vp(t),w_\vp(t))_2-(z(t),w(t))_2,\end{equation}we get, therefore, by \eqref{e3.3} and \eqref{e3.8}, 
\begin{equation}\label{e3.12}
\hspace*{-6mm}\barr{ll}
0\!\!&\le h_\vp(t)
\le h_\vp(0+){+}2\vp\!\dd\int^t_0\!\!
(\Phi_\vp(z_\vp(s)),w_\vp(s))_2ds
{-}2\alpha_3\!\dd\int^t_0
\!\!|w_\vp(s)|^2_2ds\\
&+2\alpha_1|D|_\9\dd\int^t_0\!|\nabla\Phi_\vp(z_\vp(s))|_2
|w_\vp(s)|_2
ds+2
\int^t_0|\gamma_\vp(s)|ds,\, \ff t\in[0,T].\earr\hspace*{-6mm}\end{equation} Moreover, since $z\in L^\9((0,T)\times\rrd)$ and by \eqref{e3.5} we obtain

\begin{equation}
\label{e3.13}
\vp|\Phi_\vp(z_\vp(t))|_\9\le|z_\vp(t)|_\9\le|z(t)|_\9,\mbox{\ \ a.e. }t\in(0,T).\end{equation} 
Taking into account that 
 $t\to\Phi_\vp(z_\vp(t))$ has an $L^2$ conti\-nuous version on $[0,T]$, there exists $f\in L^2$ such that 
 $$\lim_{t\to0}\Phi_\vp(z_\vp(t))=f\mbox{\ \ in }L^2.$$ 
 \n Furthermore, for every $\vf\in C^\9_0(\rrd)$, $s\in(0,T),$
 $$0\le h_\vp(s)\le|\Phi_\vp(z_\vp(s))-f|_2|z_\vp(s)|_2+|f-\vf|_2|z_\vp(s)|_2+|(\vf*\theta_\vp,z(s))_2|.$$Hence, by \eqref{e3.1},
 $$\barr{ll}
 0\le h_\vp(0+)\!\!
 &=\dd\lim_{t\downarrow0}h_\vp(t)
 =\lim_{t\to0}\
 {\rm ess}\sup_{\hspace*{-4mm}s\in(0,t)}h_\vp(s)\vsp
 &\le\left(\dd\lim_{t\to0}|\Phi_\vp(z_\vp(t))-f|_2+|f-\vf|_2\right)|z_\vp|_{L^\9(0,T;L^2)}\vsp&
 +\dd\lim_{t\to0}\ {\rm ess}\sup_{\hspace*{-4mm}s\in(0,t)}|(\vf*\theta_\vp,z(s))_2|=|f-\vf|_2|z_\vp|_{L^\9(0,T;L^2)}.\earr$$Since $C^\9_0(\rrd)$ is dense in $L^2(\rrd)$, we find 
\begin{equation}\label{c3.16'}
	 h_\vp(0+)=0.\end{equation}
 On the other hand, taking into account that, for a.e. $t\in(0,T)$,
\begin{equation}\label{e3.15}
\vp\Phi_\vp(z_\vp(t))+(-\Delta)^s\Phi_\vp(z_\vp(t))=z_\vp(t),\end{equation}we get that
\begin{equation}\label{e3.16}
\barr{r}
\vp|\Phi_\vp(z_\vp(t))|^2_2+|(-\Delta)^{\frac s2}\Phi_\vp(z_\vp(t))|^2_2=(z_\vp(t),\Phi_\vp(z_\vp(t)))_2=h_\vp(t),\vsp\mbox{ for a.e. } t\in(0,T),\earr\end{equation}and
$$
\barr{r}
\vp|(\Phi_\vp(z_\vp(t),w_\vp(t)))_2|\le\vp|\Phi_\vp(z_\vp(t))|_\9|w_\vp(t)|_1
\le|z(t)|_\9|w(t)|_1,\vsp\mbox{ for a.e. }t\in(0,T).\earr$$
We note that by   \eqref{e3.15} and Parseval's formula we have
$$|\nabla\Phi(z_\vp(t))|^2_2=\int_\rrd
\frac{|\calf(z_\vp(t))(\xi)|^2|\xi|^2}{(\vp+|\xi|^{2s})^2}\ d\xi,\ \ff t\in(0,T),$$and
$$h_\vp(t)=\int_\rrd\ \frac{|\calf(z_\vp(t))(\xi)|^2}{\vp+|\xi|^{2s}}\ d\xi,\ \ \ff t\in(0,T).$$This yields 
\begin{equation}\label{e3.15z}
\hspace*{-3mm}\barr{ll}
|\nabla\Phi_\vp(z_\vp(t))|^2_2\!\!\!
&\le
R^{2(1-s)}\dd\int_{[|\xi|\le R]}\frac{|\calf(z_\vp(t))(\xi)|^2}{\vp+|\xi|^{2s}}\ d\xi\vsp&
+\dd\int_{[|\xi|\ge R]}|\calf(z_\vp(t))(\xi)|^2|\xi|^{2(1-2s)}d\xi\vspace*{2mm}\\ 
&\le R^{2(1-s)}h_\vp(t)+R^{2(1-2s)}|z_\vp(t)|^2_2,\ \ff t\in(0,T),\ R>0,\earr\end{equation}because $2s\ge1.$ 

We shall prove now that
\begin{equation}\label{e3.15a}
\lim_{\vp\to0}\vp(\Phi_\vp(z_\vp(t)),w_\vp(t))_2=0,\mbox{ a.e. }t\in(0,T).\end{equation}
Since by \eqref{e3.13}
\begin{equation}\label{e3.18a}
\vp|(\Phi_\vp(z_\vp(t)),w_\vp(t))_2|\le|z_\vp(t)|_\9|w_\vp(t)|_1\le|z(t)|_\9|w(t)|_1,\end{equation}it suffices to show that 

\begin{equation}\label{e3.18aa}
\lim_{\vp\to0}\vp|\Phi_\vp(z_\vp(t))|_\9=0,\ \mbox{ a.e. }t\in(0,T).\end{equation}To prove \eqref{e3.18aa} we proceed similarly as in the proof of \cite[Lemma 1]{5}. 
By (A.6) in the Appendix we have for a.e. $t\in(0,T)$
$$\vp\Phi_\vp(z_\vp(t))(x)=\vp^{\frac d{2s}}\int_\rrd g^s_1(\vp^{\frac1{2s}}(x-\xi))z_\vp(t)(\xi)d\xi\mbox{ for a.e. }x\in\rrd.$$
This yields for a.e. $x\in\rrd$
$$|\vp\Phi_\vp(z_\vp(t))(x)|\le C_r\vp^{d/2s}|z(t)|_1+\vp^{d/2s}|z(t)|_\9\int_{[\vp^{1/2s}|x-\xi|\le r]}g^s_1(\vp^{1/2s}(x-\xi))d\xi,$$where $C_r:=\sup\{g^s_1(x);|x|\ge r\}\ (<\9,$ since $g^s_1\in L^\9(B_r(0)^C)$ by  (A.7)). Therefore, for a.e. $x\in\rrd$,
$$|\vp\Phi_\vp(z_\vp(t))(x)|\le C_r\vp^{d/2s}|z(t)|_1+|z(t)|_\9\int_{[|\xi|\le r]}g^s_1(\xi)d\xi.$$Since $g^s_1\in L^1$ by  (A.4), letting first $\vp\to0$ and then $r\to0$, \eqref{e3.18aa} follows, as claimed. 

By \eqref{e3.18a}, \eqref{e3.18aa} and the dominated convergence theorem, it follows that
\begin{equation}\label{e3.21}
\lim_{\vp\to0}\vp\int^t_0(\Phi_\vp(z_\vp(s)),w_\vp(s))_2ds=0,\ t\in[0,T].\end{equation}
Next, by \eqref{e3.12}, \eqref{c3.16'} and \eqref{e3.15z}, we have
$$\barr{ll}
0\le h_\vp(t)\!\!\!&\le 2\vp\dd\int^t_0|(\Phi_\vp(z_\vp(s)),w_\vp(s))_2|ds
-2\alpha_3\dd\int^t_0|w_\vp(s)|^2_2ds\vsp
&\ \ \ +2\alpha_1|D|_\9\dd\int^t_0|\nabla\Phi_\vp(z_\vp(t))|_2|w_\vp(s)|_2ds 
+2\dd\int^t_0|\gamma_\vp(s)|ds\vsp 
&\le\eta_\vp(t)+2\alpha_1|D|_\9\dd\int^t_0
\(R^{1-s}h^{\frac12}_\vp(t)+R^{1-2s}|z_\vp(t)|_2\)
|w_\vp(s)|_2ds\vsp
&\ \ \ -2\alpha_3\dd\int^t_0|w_\vp(s)|^2_2ds,\ \ff\, t\in[0,T],\ R>0,\earr$$
where
$$\eta_\vp(t):=2\vp\int^t_0|(\Phi_\vp(z_\vp(s)),w_\vp(s))_2|ds
+2\int^t_0|\gamma_\vp(s)|ds.$$This yields
$$\barr{l}
h_\vp(t)
\le\eta_\vp(t)+2\alpha_1|D|_\9
\(R^{2(1-s)}\lbb\dd\int^t_0 h_\vp(s)ds+
\dd\int^t_0\Big(R^{1-2s}|z_\vp(s)|^2_2\right.\vsp
\quad\left.\qquad
+\dd\Big(\frac1{4\lbb}+R^{1-2s}\Big)
|w_\vp(s)|^2_2\Big)ds\)
-2\alpha_3\dd\int^t_0|w_\vp(s)|^2_2ds,\ \ff\lbb>0,\ R>0.\hspace*{4mm}\earr$$Taking into account that, by \eqref{e3.3a},  
\begin{equation}\label{e3.22'}
|z(t)|_2\le \alpha^{-1}_2|w(t)|_2,\ \ \ff t\in(0,T),\end{equation}
we have
$$|z_\vp(t)|_2\le \alpha^{-1}_2|w_\vp(t)|_2+\nu_\vp(t),\ \ff t\in(0,T),$$where $\nu_\vp\to0$ in $L^2(0,T)$. Then, we get, for $\lbb$, $R>0$, suitably chosen, 
\begin{equation}\label{e3.22}
\barr{r}
\dd
0\le h_\vp(t)\le\eta_\vp(t)+C\int^t_0 h_\vp(s)ds,\ \mbox{ for } t\in[0,T],\earr\end{equation}where $C>0$ is independent of $\vp$  and $\dd\lim_{\vp\to0}\eta_\vp(t)=0$ for all $t\in[0,T]$. 

In particular, by \eqref{e3.16}, it follows that
\begin{equation}\label{e3.23}
0\le h_\vp(t)\le\eta_\vp(t)\exp(Ct),\ \ff\,t\in[0,T].\end{equation} 
This implies that $h_\vp(t)\to0$  as $\vp\to0$  for every $t\in[0,T]$, 
hence by \eqref{e3.16}   the left hand side of \eqref{e3.15} converges to zero in $S'$. Thus, $0=\lim\limits_{\vp\to0} z_\vp(t)=z(t)$ in $S'$ for a.e. $t\in(0,T)$, which implies $y_1\equiv y_2.$ If $D\equiv0$, we see by \eqref{e3.12}, \eqref{c3.16'} and \eqref{e3.15a} that $0\le h_\vp(t)\le\eta_\vp(t)$, $\ff\,t\in(0,T)$, and so the conclusion follows without invoking that $\beta'>0$, which was only used to have \eqref{e3.3a} and thus \eqref{e3.22'}.~\hfill$\Box$

\bk\n{\bf Linearized uniqueness.} In particular,   the li\-nearized uniqueness for equation \eqref{e1.10} follows  by Theorem \ref{t3.1}. More precisely,

\begin{theorem}\label{t3.2} Under assumptions of Theorem {\rm\ref{t3.1}}, let $T>0$, $u\in L^\9((0,T)\times\rrd)$  and let $y_1,y_2\in L^\9((0,T)\times\rrd)$ with $y_1-y_2\in L^1((0,T)\times\rrd)\cap L^\9(0,T;L^2)$ be two distributional solutions to the equation
	\begin{equation}
	\label{e3.30a}\barr{l}
	y_t+(-\Delta)^s\(\dd\frac{\beta(u)}uy\)+{\rm div}(yDb(u))=0\mbox{ in }\cald'((0,T)\times\rrd),\vsp
	y(0)=u_0,\earr
	\end{equation}
where $u_0$ is a measure of finite variation on $\rrd$ and $\frac{\beta(0)}0:=0$. If \eqref{e3.1} holds, then $y_1\equiv y_2$.\end{theorem}

\n{\bf Proof.} We note first that
$$\barr{c}
\dd\frac{\beta(u)}u,\ b(u)\in L^\9((0,T)\times\rrd),\vsp 
|Db(u)|_\9\le C_1\,\left|\dd\frac{\beta(u)}u\right|_\9\le C_2.\earr$$If $z=y_1-y_2,\ w=\frac{\beta(u)}u\,(y_1-y_2)$, we see that
$$\barr{ll}
wz\ge\left|\frac{\beta(u)}u\right|_\9+|w|^2,&\mbox{ a.e. on $(0,T)\times\rrd,$}\vsp
|Db(u)z|\le C_2|w|,&\mbox{ a.e. on $(0,T)\times\rrd.$}\earr$$Then, we have
$$z_t+(-\Delta)^sw+{\rm div}(Db(u)z)=0,$$
and so, arguing as in the proof of Theorem \ref{t3.1}, we get that $y_1\equiv y_2$. The details are omitted.~\hfill$\Box$


\section{Applications to McKean--Vlasov equations with L\'evy noise}\label{s4}
\setcounter{equation}{0}

\subsection{Weak existence}\label{s4.1}

To prove weak existence for \eqref{e1.11}, we use the recent results in \cite{7aa} and Theorem \ref{t2.3}.

\begin{theorem}\label{t4.1} Assume that Hypotheses {\rm(i)--(iv)} from Section {\rm\ref{s1}} hold and let $u_0\in L^1\cap L^2$. Assume that either $u_0\in\ov{D(A)}$ or that $\beta\in C^\9(\rr)$ $($see Theorem {\rm\ref{t2.3}}$)$ and let $u$ be the solution of \eqref{e1.1}  from Theorem {\rm\ref{t2.3}}.  Then, there exists a stochastic basis $\mathbb{B}:=(\ooo,\calf,(\calf_t)_{t\ge0},\mathbb{P})$ and a $d$-dimensional isotropic $2s$-stable process $L$ with L\'evy measure $\frac{dz}{|z|^{d+2s}}$ as well as an $(\calf_t)$-adapted c\`adl\`ag process $(X_t)$ on $\ooo$ such that, for
\begin{equation}\label{e4.1}
\call_{X_t}(x):=\frac{d(\mathbb{P}\circ X^{-1}_t)}{dx}\,(x),\ t\ge0,\end{equation}we have
\begin{equation}\label{e4.2}
	\barr{rcl}
dX_t&=&\dd D(X_t)b(\call_{X_t}(X_t))dx+\frac{\beta(\call_{X_t}(X_{t-}))}{\call_{X_t}(X_{t-})}\,dL_t,\\
\call_{X_0}&=&u_0.\earr\end{equation}
Furthermore, 
\begin{equation}\label{e4.2'}
	\call_{X_t}=u(t,\cdot),\ t\ge0,
	\end{equation}
in particular, $((t,x)\mapsto\call_{X_t}(x))\in L^\9([0,T]\times\rrd)$ for every $T>0.$
\end{theorem}

\n{\bf Proof.} Let $u$ be the mild (hence distributional) solution of \eqref{e1.1} from Theorem \ref{t2.3}. 

By the well known formula that
\begin{equation}\label{e4.3}
	(-\Delta)^sf(x)=-c_{d,s}{\rm P.V.}-\int_\rrd(f(x+z)-f(x))\,\frac{dz}{|z|^{d+2s}}
\end{equation}
with $c_{d,s}\in(0,\9)$ (see \cite[Section 13]{12''}), and  since, as an easy calculation shows, 
\begin{equation}
\label{e4.4'}
\barr{r}
\dd\int_A\frac{\beta(u(t,x))}{u(t,x)}\ \frac{dz}{|z|^{d+2s}}=
\dd\int_\rrd \one_A
\(\(\frac{\beta(u(t,x))}{u(t,x)}\)^{\frac1{2s}}z\)\frac{dz}{|z|^{d+2s}},\vsp A\in B(\rrd\setminus\{0\}),
\earr\end{equation}we have
\begin{equation}\label{e4.4}\barr{l}
\dd\frac{\beta(u(t,x))}{u(t,x)}\,(-\Delta)^sf(x)\\
\qquad\dd=-c_{d,s}{\rm P.V.}\int_\rrd
\(f\(x+\(\frac{\beta(u(t,x))}{u(t,x)}\)^{\frac1{2s}}z\)-f(x)\)\frac{dz}{|z|^{d+2s}}.\earr
\end{equation}
As is easily checked, Hypotheses (i)--(iv) imply that condition (1.18) in \cite{7aa} holds. Furthermore, it follows by Theorem \ref{t2.3} that
$$\mu(dx):=u(t,x)(dx),\ t\ge0,$$solves the \FP\ equation \eqref{e1.10} with $u_0:=u_0(x)dx.$ Hence, by \mbox{\cite[Theorem 1.5]{7aa},} \eqref{e4.4'}, \eqref{e4.4} and \cite[Theorem 2.26, p.~157]{15'},  there exists a stochastic basis $\mathbb{B}$ and $(X_t)_{t\ge0}$ as in the assertion of the theorem, as well as a Poisson random measure $N$ on $\rrd\times[0,\9)$ with intensity $|z|^{-d-2s}dz\,dt$ on the stochastic basis $\mathbb{B}$ such that for
\begin{equation}\label{e4.5}
L_t:=\int^t_0\int_{|z|\le1}z\wt N(dz\,ds)+\int^t_0\int_{|z|>1}z\,N(dz\,ds),\end{equation}
\eqref{e4.1}, \eqref{e4.2} and \eqref{e4.2'} hold. Here,
$$\wt N(dz\,dt):=N(dz\,dt)-|z|^{-d-2s}dz\,dt.\eqno\Box$$

\subsection{Weak uniqueness}\label{s4.2}

\begin{theorem}\label{t4.2} Assume that Hypotheses {\rm(j)--(jjj)}, resp. {\rm(j)$'$, (jj), (jjj)} if \mbox{$D\equiv0$,}  from Section {\rm\ref{s3}} hold and let $T>0$. Let $(X_t)$ and $(\wt X_t)$ be two c\`adl\`ag processes on two $($possibly different$)$ stochastic bases $\mathbb{B}, \wt{\mathbb{B}}$ that are weak solutions to \eqref{e4.2} with $($possibly different$)$ $L$ and $\wt L$ defined as in \eqref{e4.5}. Assume that 
	\begin{equation}\label{e4.6}
		\((t,x)\mapsto \call_{X_t}(x)\),
		\((t,x)\mapsto\call_{\wt X_t}(x)\)\in L^\9((0,T)\times\rrd). \end{equation}Then $X$ and $\wt X$ have the same laws, i.e.,
	$$\mathbb{P}\circ X\1=\wt{\mathbb{P}}\circ\wt X\1.$$  
\end{theorem}

\n{\bf Proof.} Clearly, by Dynkin's formula both
$$\mu_t(dx):=\call_{X_t}(x)dx\mbox{\ \ and\ \ }
\wt\mu_t(dx):=\call_{\wt X_t}(x)dx$$solve the \FP\ equation \eqref{e1.10} with the same initial condition $u_0(dx):=u_0(x)dx$, hence satisfy \eqref{e3.1} with $y_1(t):=\call_{X_t}$ and $y_2(t):=\call_{\wt X_t}$. Hence, by Theorem \ref{t3.1},
$$\call_{X_t}=\call_{\wt X_t}\mbox{\ \ for all }t\ge0,$$since $t\mapsto\call_{X_t}(x)dx$ and $t\mapsto\call_{\wt X_t}(x)dx$ are both narrowly  continuous and are probability measures for all $t\ge0$, so both are in $L^\9(0,T;L^1\cap L^\9)\subset L^\9(0,T;L^2)$. 

Now, consider the linear \FP\ equation
\begin{equation}\label{e4.7}
\barr{c}
v_t+(-\Delta)^s\beta(\call_{X_t})+{\rm div}(Db(\call_{X_t})v)=0,\vsp
v(0,x)=u_0(x),\earr\end{equation}
again in the weak (distributional) sense analogous to \eqref{e1.10}. Then, by Theo\-rem \ref{t3.2} we conclude that $\call_{X_t}$, $t\in[0,T]$, is the unique solution to \eqref{e4.7} in  $L^\9(0,T;L^1\cap L^\9)$. Again by Dynkin's formula, both $\mathbb{P}\circ X^{-1}$ and  $\wt{\mathbb{P}}\circ \wt X\1$ solve the martingale problem with initial condition $u_0(dx):=u_0(x)dx$ for the linear Kolmogorov operator
\begin{equation}\label{e4.8}
K_{\call_{X_t}}:=-\frac{\beta(\call_{X_t})}{\call_{X_t}}\,(-\Delta)^s+b(\call_{X_t})D\cdot\nabla.\end{equation}
Since the above is true for all $u_0\in L^1\cap L^\9$, and also holds when we consider \eqref{e1.1} resp \eqref{e1.10} with start in $s>0$ instead of zero, it follows by exactly the same arguments as in the proof of Lemma 2.12 in \cite{13'}     that
$$\mathbb{P}\circ X\1=\wt{\mathbb{P}}\circ\wt X\1.\eqno\Box$$

\begin{remark}\label{r4.3}\rm Let for $s\in[0,\9)$ and $\calz:=\{\zeta\equiv\zeta(x)dx\mid\zeta\in L^1\cap L^\9,\ \zeta\ge0,$ $|\zeta|_1=1\}$
	$$\mathbb{P}_{(s,\zeta)}:=\mathbb{P}\circ X^{-1}(s,\zeta),$$where $(X_t(s,\zeta))_{t\ge0}$ on a stochastic basis $\mathbb{B}$ denotes the solution of \eqref{e1.11} with initial condition $\zeta$ at $s$. Then, by Theorems \ref{t3.1}, \ref{t3.2} and \ref{t4.2}, exactly the same way as Corollary 4.6 in \cite{11''}, one proves that $\mathbb{P}_{(s,\zeta)},$ $(s,\zeta)\in[0,\9)\times\calz,$ form a nonlinear Markov process in the sense of McKean (see \cite{9'}).
	\end{remark}

\begin{remark}\label{r4.4} \rm \eqref{e4.2'} in Theorem \ref{t4.1} says that our solution $u$ of \eqref{e1.1} from Theorem \ref{t2.3} is the law density of a c\`adl\`ag process solving \eqref{e4.2} or resp. by Remark \ref{r4.3} above that it is the law density of a nonlinear Markov process. This realizes McKean's vision formulated in \cite{9'} for solutions to nonlinear parabolic PDE. So, our results show that it is also possibly for nonlocal PDE of type \eqref{e1.1}.\end{remark}

\begin{remark}\label{r4.5}\rm In a forthcoming paper \cite{6'}, we achieve similar results as in this paper in the case where $(-\Delta)^s$ is replaced by $\psi(-\Delta)$, where $\psi$ is a Bernstein function (see \cite{12''}). 
	\end{remark}


\subsection*{Appendix:\\ Representation and properties of the integral kernel of $(\vp I+(-\Delta)^s)\1$}

Let $s\in(0,1)$ and let $\calf:S'(\rrd)\to S'(\rrd)$ be the Fourier transform, as defined in \eqref{e1.3}. 

It is well known (see, e.g., \cite[Chap.~II, Sect.~4c]{MR92}) that for $t>0$ the integral kernel $p^s_t$ of the operator $T^s_t:=\exp(-t(-\Delta)^s)$ is related to the kernel $p_t$ of $\exp(t\Delta)$ by the following subordination formula
$$p^s_t(x)=\int^\9_0p_r(x)\eta^s_{t}(dr),\ x\in\rrd,\eqno{\rm(A.1)}$$ 
where $(\eta^s_t)_{t>0}$ is the one-sided stable semigroup of order $s\in(0,1)$, which is defined through its Laplace transform by

$$\int^\9_0 e^{-\lbb r}\eta^s_t(dr)=e^{-t\lbb^s},\ \lbb>0.\eqno{\rm(A.2)}$$ 
Furthermore, since $(\vp I+(-\Delta)^s)\1$ is the Laplace transform of the semigroup $T^{(s)}_t,$ $t\ge0$, it follows that
$$g^{s}_\vp(x)=\int^\9_0e^{-\vp t}\int^\9_0p_r(x)\eta^s_{t}(dr)dt,\ x\in\rrd,\eqno{\rm(A.3)}$$ 
is the integral kernel of $(\vp I+(-\Delta)^s)\1$. Obviously, since each $\eta^s_t$ is a probability measure, we have
$$\vp\int_\rrd g^s_\vp(x)dx=1.\eqno{\rm(A.4)}$$
Since
$$p_r(x)=\dd\frac1{(4\pi r)^{d/2}}\ e^{-\frac 1{4r}\,|x|^2}
=(2\pi)^{-d}\dd\int_\rrd e^{i\<x,y\>}e^{-r|y|^2}dy,\ x\in\rrd,\eqno{\rm(A.5)}$$
plugging the first equality in (A.5) into (A.3),  we obtain

$$g^s_\vp(x)=\int^\9_0e^{-\vp t}\int^\9_0
\frac1{(4\pi r)^{d/2}}\ e^{-\frac1{4r}\,|x|^2}
\eta^s_t(dr)dt.$$It is well known and trivial to check from the definition that for $\gamma\in(0,\9)$ the image measure of $\eta^s_t$ under the map $r\mapsto \gamma r$ is equal to $\eta^s_{\gamma^st}$. Hence, by an elementary computation we find
$$g^s_\vp(x)=|x|^{-d+2s}\int^\9_0\int^\9_0 
e^{-\vp|x|^{2s}t}
\frac1{(4\pi rt^{1/s})^{d/2}}\  e^{-\frac1{4rt^{1/s}}}\,dt\,\eta^s_1(dr),\eqno{\rm(A.6)}$$which, since $s<\frac d2$ and $\eta^s_1$ is a probability measure, in turn implies
$$g^s_\vp(x)\sim|x|^{-d+2s}\mbox{\ \ as }|x|\to0$$and
$$g^s_\vp\in L^\9(\rrd\setminus B_R(0)),\ \ff R>0.\eqno{\rm(A.7)}$$
Plugging the second equality in (A.5) into (A.3), it follows by (A.2) that
$$\barr{ll}
g^s_\vp(x)\!\!\!
&=(2\pi)^{-d}\dd\int^\9_0
e^{-\vp t}\int_\rrd e^{i\<x,y\>}e^{-t|y|^{2s}}\,dy\,dt\vsp
&=(2\pi)^{-d/2}\dd\int^\9_0
\calf\(e^{-t(\vp+|\,\cdot\,|^{2s})}\)(x)dt,\ x\in\rrd.\earr\eqno{\rm(A.8)}$$Hence  
$$g^s_\vp=(2\pi)^{-d/2}\calf\(\frac1{\vp+|\cdot|^{2s}}\)\mbox{\ \ in }S'(\rrd). \eqno{\rm(A.9)}$$ 
Finally, from (A.8) it follows that
$$g^s_\vp(x)=\vp^{\frac{d-2s}{2s}}g^s_1
\(\vp^{\frac1{2s}}x\),\ x\in\rrd. \eqno{\rm(A.10)}$$ 
\n{\bf Acknowledgement.} This work was supported by the DFG through SFB 1283/2 2021-317210226 and by a grant of the Ministry of Research, Innovation and Digitization, CNCS--UEFISCDI project  PN-III-P4-PCE-2021-0006, within PNCDI III. A part of this work was done during a very pleasant stay of the second named author at the University  of Madeira as a guest of Jos\'e Luis da Silva. We are grateful for his hospitality and for many discussions as well as for carefully reading large parts of this paper.

\end{document}